\documentclass[12pt]{article}

\usepackage[utf8]{inputenc}
\usepackage[english]{babel}
\usepackage{graphicx}
\usepackage{amsmath, accents}
\usepackage{amssymb}
\usepackage{amsfonts}
\usepackage{amsthm}
\usepackage{mathrsfs}
\usepackage{latexsym}
\usepackage{subfigure}
\usepackage{color}
\usepackage{xcolor}
\usepackage{url,hyperref}
\usepackage[noend]{algpseudocode}

\newcommand{\Z}{\mathbb{Z}}
\newcommand{\R}{\mathbb{R}}
\newcommand{\abs}[1]{\left\vert #1\right\vert}

\newtheorem{theorem}{Theorem}
\newtheorem{ex}[theorem]{Example}
\newtheorem{defin}[theorem]{Definition}
\newtheorem{prop}[theorem]{Lemma}

\newtheorem{cor}[theorem]{Corollary}

\begin{document}

\title{A linear time approach to three-dimensional reconstruction by discrete tomography}

\author{Matthew Ceko\thanks{School of Physics and Astronomy, Monash University, Melbourne, Australia} \and Silvia M.C.~Pagani\thanks{Dipartimento di Matematica e Fisica ``N.~Tartaglia'', Università Cattolica del Sacro Cuore, Brescia, Italy. Corresponding author; email: silvia.pagani@unicatt.it} \and Rob Tijdeman\thanks{Mathematical Institute, Leiden University, Leiden, The Netherlands}}

%

\maketitle

\begin{abstract}
The goal of discrete tomography is to reconstruct an unknown function $f$ via a given set of line sums. In addition to requiring accurate reconstructions, it is favourable to be able to perform the task in a timely manner. This is complicated by the presence of ghosts, which allow many solutions to exist in general. In this paper we consider the case of a function $f : A \to \mathbb{R}$ where $A$ is a finite grid in $\mathbb{Z}^3$. Previous work has shown that in the two-dimensional case it is possible to determine all solutions in parameterized form in linear time (with respect to the number of directions and the grid size) regardless of whether the solution is unique. In this work, we show that a similar linear method exists in three dimensions under the condition of nonproportionality.
We show that the condition of nonproportionality is fulfilled in the case of three-dimensional boundary ghosts.\smallskip

\textbf{Keywords}: discrete tomography; ghost; lattice direction; linear time algorithm; three-di\-men\-sion\-al re\-con\-struc\-tion.
\end{abstract}

\section{Introduction}
Tomography is the process of reconstructing an object from a set of its projected views. The usual continuous tomography requires many projections for accurate solutions. Imaging sensitive objects such as fine nanostructures or biological matter often limit the number of projections that can be taken, as to not perturb or destroy the object. This is the case as well for objects which change quickly in time. The methods of discrete tomography permit useful reconstructions or approximations to be obtained with relatively little projection information. Discrete tomography considers a function $f$, representing the object, defined on a finite grid $A$ of $\mathbb{Z}^2$. Instead of continuous integral projections, discrete line sums are used. These sum the $f$-values at grid points along lines in a small number of directions. In this work, we assume that these line sums are free of noise and errors. Algorithms based on this condition have found use in image processing and data security, cf.~\cite{csgkn, gue, per, vbrev}.

Discrete tomography began from analysing the problem of recovering binary matrices from their row and column sums. In 1957, Ryser gave a necessary and sufficient condition, and an algorithm for reconstruction of this problem, \cite{rys}. In 1978, Katz showed that reconstruction is possible for any number of directions in the absence of a nontrivial function with vanishing line sums over the set of directions, known as a ghost (or switching function), \cite{katz}. Hajdu and Tijdeman offered an algebraic interpretation of discrete tomography in 2001, which viewed ghosts as polynomials, \cite{hatij}. Ghosts of minimal size are termed minimal ghosts. In \cite{hatij} it was shown that every ghost is a linear combination of minimal ghosts. Hence, points of the grid $A$ have uniquely determinable $f$-values if they are not in the union of all the ghost domains. Furthermore, assigning arbitrary function values to a suitable subset of the ghost domain induces unique function values for all other points. This fact was exploited by Dulio and Pagani in \cite{dupa}, wherein a rounding theorem was proven which allowed exact and unique binary tomographic reconstructions from the minimum Euclidean norm solution. It also motivated the construction of boundary ghosts by Ceko, Petersen, Svalbe and Tijdeman, \cite{cpst}. Boundary ghosts are ghosts with a thin annulus as domain, having a large interior of points where the $f$-values are uniquely determined. Ceko and Tijdeman \cite{ct} showed that a straightforward generalization to three dimensions does not exist, but that a similar structure exists by combining three recursions. In this paper we introduce a reconstruction method for such three-dimensional boundary ghosts.

An important aspect of discrete tomography is the time in which reconstructions can be obtained. It was shown by Gardner, Gritzmann and Prangenberg in 1999 that a function $f: A \to \mathbb{N}$ can be reconstructed in linear time for two directions, but that the problem is NP-complete for $s \geq 3$ directions, \cite{ggp99}. They also showed that the problem is NP-complete for $s \geq 2$ when the $f$ attains six or more different values, \cite{ggp00}. The situation is completely different if the co-domain is an integral domain. In 2015 Dulio, Frosini and Pagani proved that the function values in the corners of $A$ can be uniquely determined in linear time for $s=2$ in \cite{deda}, and gave conditional results for $s=3$ in \cite{dgci2016,3dirext}. This result was generalised by Pagani and Tijdeman for all $s$ in \cite{pati}. They showed that the function values at $A$ outside of the convex hull of the union of all the ghost domains can be computed in linear time. If there is no ghost, all function values may be obtained in linear time.

The coalescence of these results shows that an algorithm can be composed, which allows for reconstruction of all solutions from the line sums in linear time in terms of the product of the number of directions and the size of the grid. Such an algorithm was provided in \cite{cst}. In the present work, we extend the algorithm to three dimensions under the condition of nonproportionality.
This is achieved through repeated applications of the aforementioned 2D algorithm along coordinate planes in 3D. We show that three-dimensional boundary ghosts have a nonproportional set of primitive directions so that a 3D discrete reconstruction can be computed in linear time with respect to the product of the grid size and number of directions. We illustrate this with an example of 11 directions.

In the next section we introduce notation and the main result. The relevant results of the two-dimensional case are given in Section \ref{sec:2dim}. Here we prove that in the algorithm of \cite{cst} primitivity may be replaced by nonproportionality. In Section \ref{nonp} we show that the three-dimensional boundary ghosts satisfy the nonproportionality property. We treat the three-dimensional case in Section \ref{dimthree} in four steps where we illustrate each step by the example of the boundary ghost. Section \ref{sec:complexity} offers an analysis of the computational complexity of the algorithm. Finally we state some conclusions in the last section.

\section{Notation and main result} \label{dim3}

Throughout the paper we consider an $X$ by $Y$ by $Z$ grid of points
$$A = \{(x,y,z) \in \mathbb{Z}^3 : 0 \leq x < X, 0 \leq y < Y, 0 \leq z <Z\}.$$

\begin{defin}
A triple $(a,b,c) \in \mathbb{Z}^3$ with $(a,b,c) \neq (0,0,0)$ is called a direction, where we identify $(a,b,c)$ and $(-a,-b,-c)$.
It is called primitive if $\gcd(a,b,c)=1$.
\end{defin}

Let $D$ be a set of directions $d_h = (a_h,b_h,c_h)$ for $h = 1,2, \dots, s$.

\begin{defin} \label{np}
The set $D$ is called nonproportional if, for $h=1,\ldots s$,
\begin{itemize}
\item[$1)$] the ratios $a_h:c_h$ with $a_hc_h \not= 0$ are distinct,
\item[$2)$] the ratios $b_h:c_h$ with $c_h \not= 0$ are distinct,
\item[$3)$] the ratios $a_h:b_h$ with $c_h=0$ are distinct.
\end{itemize}
\end{defin}

\begin{defin}
A lattice line $L$ is a set of points $L=L(a,b,c, x_0,y_0,z_0) = \{(x_0,y_0,z_0) + t(a,b,c) : t \in \mathbb{Z}\}$ where $(x_0,y_0,z_0) \in A$ and $(a,b,c)$ is a direction.

For a function $f : A \to \mathbb{R}$, the line sum of $f$ along $L$ is defined as
$$\ell(a,b,c,x_0,y_0,z_0,f) = \sum_{(x,y,z) \in L \cap A} f(x,y,z).$$
\end{defin}

If $\gcd(a,b,c)>1$, the lattice line does not contain all integer points on the line $\{(x_0,y_0,z_0) + t(a,b,c) : t \in \mathbb{R}\}$.

\begin{defin}
A nontrivial function $g: A \to \mathbb{R}$ is called a ghost of $(A,D)$ if all the line sums of $g$ in all the directions of $D$ are $0$. The support of a ghost is called its ghost domain. A ghost of minimal size is called minimal ghost.
\end{defin}

\begin{defin}
We call $(A,D)$ valid if $\sum_{h=1}^s \abs{a_h} < X$, $ \sum_{h=1}^s \abs {b_h}  < Y$ and $\sum_{h=1}^s \abs{c_h} < Z$, and nonvalid otherwise.
\end{defin}

\begin{defin}
An elementary operation is an addition, subtraction, multiplication, division, a determination of the largest of two given quantities, or an assignment. We call an algorithm linear if the number of required elementary operations is $\mathcal{O}(sXYZ)$, where $s,X,Y,Z$ are defined as above.
\end{defin}

Let the set $D$ consist of primitive directions. The function $f: A \to \mathbb{R}$ is uniquely determined by its line sums in the directions of $D$ if and only if $(A,D)$ is nonvalid, \cite{ht07, katz}. If $(A,D)$ is valid, then there is a minimal ghost $ g_{(0,0,0)} : A_{(0,0,0)} \to \mathbb{R}$ with
$$A_{(0,0,0)} \subseteq  \left[0,\sum_h \abs{a_h}\right] \times \left[0,\sum_h \abs{b_h}\right] \times \left[0,\sum_h \abs{c_h}\right] \cap \mathbb{Z}^3.$$
This minimal ghost is unique apart from a multiplicative factor. In fact,
\begin{equation} \label{elemt}
A_{(0,0,0)} =\left(\sum_{a_h<0} |a_h|, \sum_{b_h<0} |b_h|, \sum_{c_h<0} |c_h|\right) + \varepsilon_1 d_1 +  \varepsilon_2  d_2 + \ldots + \varepsilon_s d_s
\end{equation}
with $\varepsilon_1, \varepsilon_2,  \dots, \varepsilon_s \in \{0,1\}$.
Let $g_{(x,y,z)}$ be the corresponding minimal ghost with shifted domain $A_{(0,0,0)} + (x,y,z)$. Then every ghost is of the form $\sum_{(x,y,z) \in U} c_{(x,y,z)}g_{(x,y,z)}$ with $c_{(x,y,z)} \in \mathbb{R}$ and $U$ the subset of $(x,y,z) \in A$ for which the domain of $g_{(x,y,z)}$ fits into $A$, \cite{ht07}. It follows that if $(A,D)$ is valid, there are
\begin{equation} \label{dimspace}
\left(X- \sum_{h=1}^s \abs{a_h} \right) \times \left(Y - \sum_{h=1}^s |b_h| \right) \times \left(Z - \sum_{h=1}^s |c_h| \right)
\end{equation}
linearly independent functions $f^*: A \to \mathbb{R}$ with the same line sums as $f$ has. In so many points of $A$ the value of $f^*$ can be chosen in $\mathbb{R}$. A function $f$ is uniquely determined by its line sums in the directions of $D$ outside the union $T$ of its ghost domains. It follows from \eqref{elemt} and the definition of $U$ that
\begin{equation} \label{ghostT}
\displaystyle T = \bigcup_{i=0}^{X-1- \sum_h |a_h|}\,\, \bigcup_{j=0}^{Y- 1-\sum_h |b_h|}\,\, \bigcup_{k=0}^{Z-1 - \sum_h |c_h|} A_{(i,j,k)}.
\end{equation}

The aim of the present paper is to prove that the following algorithm can be composed with $\mathcal{O}(sXYZ)$ elementary operations.\smallskip

\noindent {\bf Algorithm A}
\begin{algorithmic}
    \Require{A set $A = \{(x,y,z) \in \mathbb{Z}^3 : 0 \leq x <X, 0 \leq y<Y, 0 \leq z <Z \}$, a nonproportional set of directions $D$ and all the line sums in the directions of $D$ of a function $f : A \to \R.$}
    \Ensure{A function $f^* : A \to \R$ which satisfies the line sums of $f$.}
\end{algorithmic}

As shown in \cite{ht07} this implies a parametric characterisation of all the functions of which the line sums in the directions of $D$ agree with those of $f$. As an application we show that the three-dimensional boundary ghosts described in \cite{ct} satisfy the conditions of the algorithm and show how the algorithm works for such functions.

\section{The two-dimensional case}\label{sec:2dim}
Let $A$ be an $X$ by $Y$ grid of points
$$A = \{(x,y) \in \mathbb{Z}^2 : 0 \leq x < X, 0 \leq y < Y\}.$$
We use the notation of Section \ref{dim3}, but since the third component of the points is 0, we omit it throughout the present section. For two dimensions nonproportional means that condition 3) of Definition \ref{np} is satisfied, that is, $ab' \not= a'b$ for any two directions $(a,b)$ and $(a',b')$. Thus, in the two-dimensional case, primitive implies nonproportional.

Let $D$ be a set of primitive directions. The function $f: A \to \mathbb{R}$ is uniquely determined by its line sums in the directions of $D$ if and only if $(A,D)$ is nonvalid, \cite{katz}. If $(A,D)$ is valid, then there are
\begin{equation} \label{numbers2}
\left( X- \sum_{h=1}^s |a_h| \right)  \left( Y- \sum_{h=1}^s |b_h| \right)
\end{equation}
linearly independent functions $f^*: A \to \mathbb{R}$ with the same line sums as $f$ has, that is, at so many points of $A$ the value of $f^*$ can be chosen in $\mathbb{R}$.  A function $f$ is uniquely determined by its line sums in the directions of $D$ outside the union $T$ of its ghost domains and not elsewhere (cf.~\cite{ht07}). We denote the convex hull of $T$ by $C$. In \cite{cst} the following algorithm is presented.
\vskip.2cm
\noindent {\bf Algorithm B} (Algorithm of \cite{cst}, Section 6.)
\begin{algorithmic}
    \Require{A set $A = \{(x,y) \in \mathbb{Z}^2 : 0 \leq x <X, 0 \leq y<Y \}$, a finite set of primitive directions $D$ and all the line sums in the directions of $D$ of a function $f : A \to \R.$}
    \Ensure{A function $f^* : A \to \R$ which satisfies the line sums.}
\end{algorithmic}
\vskip.2cm

\noindent Algorithm B computes the $f$-values of the integer points outside $T$. It uses $\mathcal{O}(sXY)$ elementary operations to compute $f^*$. The $f^*$-values themselves are calculated by subtractions only.

In the proof that Algorithm A can be made, Algorithm B is applied to the projection of $C$ along the $x$-axis, the $y$-axis and the $z$-axis. The notion `nonproportional' is introduced to guarantee that in all three applications we have nonproportional directions. As an example, the primitive direction $(5,-5,4)$ yields, when projected along the $z$-axis, the nonprimitive direction $(5,-5)$. Since we do not want to exclude directions like $(5,-5,4)$ in Algorithm A, we prove that the condition `primitive' in Algorithm B can be relaxed to `nonproportional'.

\begin{theorem} \label{nonprop}
In Algorithm B `primitive' can be replaced by `nonproportional'.
\end{theorem}

In order to show why this extension of Algorithm B is possible we restrict our attention to one corner of $A$. The same argument may be adapted to the other corners. Let $D$ be a set of distinct directions $(a_1, b_1), \ldots, (a_k,b_k)$ with $k \geq 2$, where $a_1, \ldots, a_k$ are positive integers and $b_1, \dots, b_k$ negative integers ordered such that
\begin{equation}\label{ratios}
\frac{\abs{b_1}}{a_1}<\frac{\abs{b_2}}{a_2} < \ldots < \frac {\abs{b_k}}{a_k},
\end{equation}
where $1/0$ is considered as $\infty$.

\begin{defin} \label{def:order2}
We call the vertices of $C$ $$\left(\sum_{h=1}^k a_h,0\right) , \,  \left(\sum_{h=2}^k a_h, |b_1| \right),\, \ldots,\, \left(0,\, \sum_{h=1}^k |b_h|\right) $$
the border points $(P_0,Q_0), (P_1,Q_1) \ldots, (P_k,Q_k)$, respectively.
\end{defin}

\begin{defin} \label{def:weight}
For a point $(x,y) \in A$ we define its weight $w(x,y)$ by
$$w(x,y) = \frac {a_Hy+|b_H|x}{a_H(\sum_{h=1}^H |b_h|) + |b_H|(\sum_{h=H+1}^k a_h)},$$
where $H$ is determined by
\begin{equation*}
\frac{\sum_{h=1}^{H-1} |b_h|}{\sum_{h=H}^{k} a_h} \leq \frac yx \leq \frac{\sum_{h=1}^H |b_h|}{\sum_{h=H+1}^{k} a_h}.
\end{equation*}
\end{defin}

\begin{figure}
\centering
\includegraphics[scale=0.6]{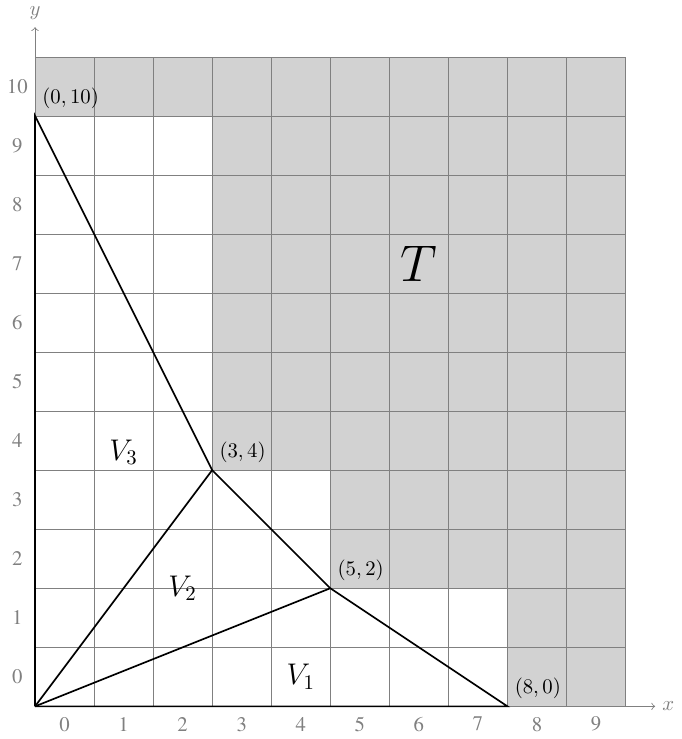}
\caption{\small The triangles $V_1,V_2,V_3$ for the set $D=\{(3,-2),(2,-2),(3,-6)\}$. The border points are $(P_0,Q_0)=(8,0)$, $(P_1,Q_1)=(5,2)$, $(P_2,Q_2)=(3,4)$, $(P_3,Q_3)=(0,10)$. They are the vertices of the convex hull $C$ of $T$ (elements of $T$ are coloured in the figure). For every $H$ the line through $(P_{H-1},Q_{H-1})$ and $(P_H,Q_H)$ is an edge of triangle $V_H$, and intersects each other triangle $V_h$, since $C$ is convex.}
\label{fig:triangles1}
\end{figure}

Obviously $w(x,y)=1$ if $(x,y)$ is one of the border points. Observe that the function $w$ is linear on $V_H$ and that the lines with equation
\begin{displaymath}
y\sum_{h=H+1}^k a_h = x \sum_{h=1}^{H} |b_h| \qquad \text{for }H=0,1,\ldots,k
\end{displaymath}
connect the origin with the border points. The lines generate $k$ full triangles $V_1,V_2,\ldots,V_k$ with as vertices of $V_H$ the origin and two consecutive border points $(P_{H-1},Q_{H-1})$ and $(P_H,Q_H)$ for $H=1,2, \dots, k$ (see Figure \ref{fig:triangles1}). For a point $(x,y)\in V_1\cup\ldots\cup V_k$ the index $H$ appearing in Definition \ref{def:weight} is that of triangle $V_H$ containing $(x,y)$.

Denote by $\delta((x,y),(x',y'))$ the distance between the points $(x,y)$ and $(x',y')$. The weight has the following property.

\begin{prop} \label{lem:weight}
For $H=1,2, \dots, s$ and $(x,y) \in V_H$ denote the intersection of the line through the origin and $(x,y)$ and the line through $(P_{H-1},Q_{H-1})$ and $(P_H,Q_H)$ by $\overline{(x,y)}$.
Then
\begin{equation*}
w(x,y) = \frac {\delta((0,0), (x,y))} {\delta((0,0),\overline{(x,y)})}.
\end{equation*}
\end{prop}

\begin{proof} We have $w(P_{H-1},Q_{H-1}) = w(P_H,Q_H) = 1$. The weight function is linear on $V_H$, hence equal to 1 on the segment connecting $(P_{H-1},Q_{H-1})$ and $(P_H,Q_H)$. If $(x,y)$ is on the line through $(0,0)$ and $(P_H,Q_H)$, then $\overline{(x,y)} = (P_H,Q_H)$ and the statement follows from the linearity of the weight function. Since on $V_H$ the weight function $w(x,y)$ is constant on lines parallel to the line connecting $(P_{H-1},Q_{H-1})$ and $(P_H,Q_H)$, the statement holds for all points in $V_H$.
\end{proof}

\begin{figure}
\centering
\includegraphics[scale=0.7]{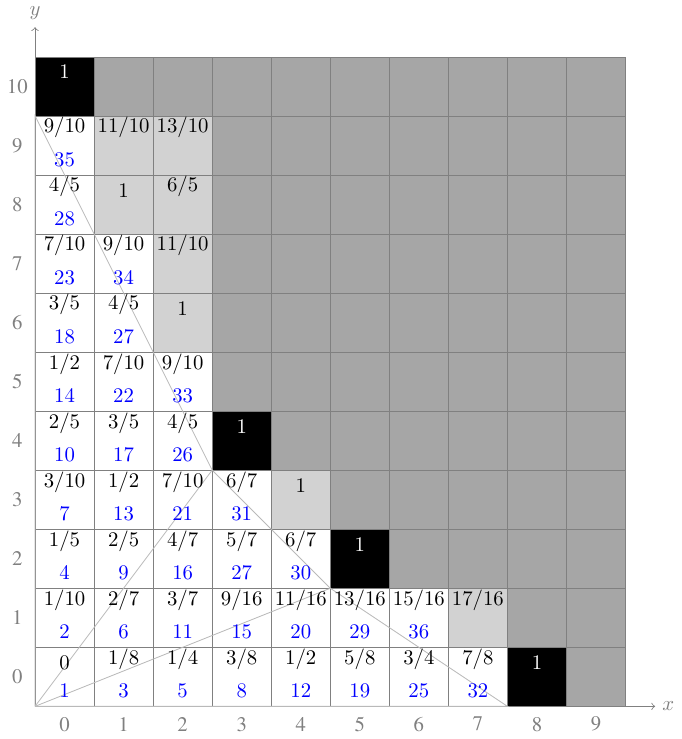}
\caption{\small The weights (upper numbers inside each pixel) and order numbers (lower numbers) for directions as in Figure \ref{fig:triangles1}. The (black) border points are in $T$ and have weight $1$. The dark points belong to $T$ as well. The points to the left of the line segments connecting consecutive border points have weights less than 1, those to the right have weights greater than 1. Since $D$ has nonprimitive directions, there are points between two consecutive border points with weight $1$. Light coloured points belong to $C$, but not to $T$. \label{fig:weight}}
\end{figure}

The reconstruction method orders the integer points of $A$ outside $C$ according to increasing weight, see Figure \ref{fig:weight}. We shall show that if such a point $(x,y)$ is in $V_H$, then it is the point with highest weight among the points on the lattice line through $(x,y)$ in direction $(a_H,b_H)$. Since the $f$-values of all other points on this lattice line are already known, the $f$-value of $(x,y)$ can be computed by subtracting all the known $f$-values of other points on that lattice line from the corresponding line sum. This is expressed in the following result.

\begin{prop} \label{bROU}
Let $A$ be the $X \times Y$ grid of points $A=\{(x,y)\in\Z^2:0 \leq x< X,\,0 \leq y< Y\}$. Let $D$ be a set of nonproportional directions $\{(a_1, b_1), \ldots, (a_k,b_k)\}$ where $a_1, \ldots, a_k$ are positive integers and $b_1, \ldots b_k$ negative integers ordered as in \eqref{ratios}. Then $(x,y) \in V_H$, not a border point, has a larger weight than all the points $(x,y) + t(a_H,b_H) \in A$ with $t \in \mathbb{Z}, t \not= 0$.
\end{prop}

This lemma for primitive directions is in \cite{pati}. We give a simpler proof for Lemma \ref{bROU} based on Lemma \ref{lem:weight}.

\begin{proof}
Consider a point $(x,y) \in V_H$ and a point $(x',y') = (x,y) +t(a_H,b_H) \in A$ with $t \in \mathbb{Z}, t \not= 0$. The triangle $V_H$ has an edge with vertices $(P_{H-1},Q_{H-1})$ and $(P_H,Q_H)$. Since $(P_{H-1},Q_{H-1})-(P_H,Q_H)= (a_H,b_H)$, the shift of $(x,y)$ by a nonzero multiple of $(a_H,b_H)$ cannot be in $V_H$. This implies that $(x',y') \in V_h$ for some $h \not= H$.
By the convexity of $C$ we have $$w(x',y') = \frac{\delta((0,0),(x',y'))}{\delta((0,0), \overline{(x',y')})} <  \frac{\delta((0,0),(x,y))}{\delta((0,0), \overline{(x,y)})} = w(x,y).$$
\end{proof}

\begin{figure}
\centering
\includegraphics[width=\columnwidth]{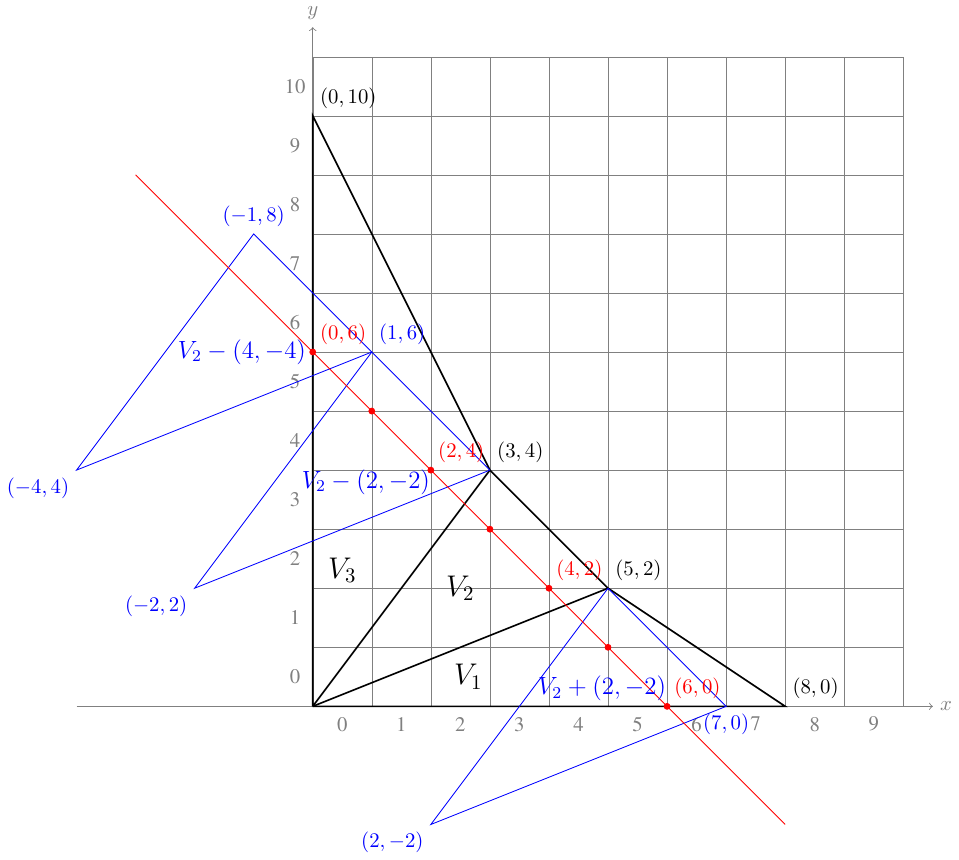}
\caption{\small The triangles $V_1,V_2,V_3$ for the set $D$ as in Figures 1 and 2. The translates $V_2\pm(2,-2)$ have only border points in common with $V_2$. The points $(4,2)+(2,-2)$, $(4,2)-(2,-2)$, $(4,2)-2(2,-2)$ have lower weights (.750, .800, .600, respectively) than $(4,2)$ has (.857). Note that some lattice points on the line $x+y=6$ are not included in Lemma \ref{bROU}. For example, point $(3,3)$ has also weight .857.}
\label{fig:triangles3}
\end{figure}

\begin{proof}[Proof of Theorem \ref{nonprop}] Replace Lemma 3 in the proof of Theorem 7 of \cite{pati} by Lemma \ref{bROU} of the present paper (cf.~Figure \ref{fig:triangles3}).
\end{proof}

\section{An infinite set of nonproportional directions} \label{nonp}
We show that the directions of the boundary ghosts from \cite{ct} are nonproportional. In the next section we use such a ghost to illustrate the method.

Define a sequence of directions by $v_0=(0,1,0), v_1=(1,0,0), v_2 = (1,1,0), v_3 = (1,0,1)$ and further, for $n=1,2, \dots$,
\begin{align*}
    v_{3n+1} &= v_{3n-1} - 2v_{3n-2}, \\
    v_{3n+2} &= v_{3n+1} - 2v_{3n-1}, \\
    v_{3n+3} &= v_{3n+1} - 2v_{3n}.
\end{align*}

A ghost generated by $(v_n)_{n=0}^N$ for some $N$ is called a boundary ghost. The theory in this paper is independent of the restriction to an initial subset and holds for any subset of $(v_n)_{n=0}^{\infty}$.

Let $v_n = (a_n,b_n,c_n)$ for $n = 1,2, \dots$. By induction it follows that $a_n$ and $b_n$ are odd for $n>3$ and $c_n=0$ if $n\not\equiv0\mod3$.

\begin{prop}
The directions $(v_n)_{n=0}^{\infty}$ are nonproportional.
\end{prop}

\begin{proof}
We distinguish between $n\equiv0\pmod3$ and $n\not\equiv0\pmod3$.

If $n=3m$, then $a_{3m}$ and $b_{3m}$ are odd for all $m>1$ and $|c_{3m}|=2^{m-1}$ for $m\geq1$. It follows that all ratios $a_{3m} : c_{3m}$ and all ratios $b_{3m} : c_{3m}$ are distinct, so conditions 1) and 2) of Definition \ref{np} are satisfied (in this case condition 3) does not occur).

Let $n\not\equiv0\pmod3$. Then $c_n=0$. Suppose $a_{h_1}:b_{h_1} = a_{h_2}:b_{h_2}$ for some $h_1>h_2$ with $h_1,h_2\not\equiv0\pmod3$. Define $(w_n)_{n=1}^{\infty}$, with $w_n=(a'_n,b'_n,0)$, in the following way:
\begin{displaymath}
w_{2m+1} = v_{3m+1}, \qquad w_{2m+2} = v_{3m+2} \qquad \text{for }m=0,1,\ldots.
\end{displaymath}
Then $w_{n+1} = w_n -2w_{n-1}$ for $n \geq 2$ and the recursion of \cite{cpst} can be applied. Formulas (5)-(6) of \cite{cpst} state that
$$a'_{n-1} = b'_n = \frac{1}{\sqrt{-7}}(\alpha^n - \overline{\alpha}^n), \text{ with } \alpha = \frac12 + \frac12 \sqrt{-7} $$
for all $n>0$. Thus $a_{h_1}:b_{h_1} = a_{h_2}:b_{h_2}$ implies $a'_kb'_{\ell} = a'_{\ell}b'_k$ for some $k>\ell$, which in turn gives
$$(\alpha^{k+1}-\overline{\alpha}^{k+1}) (\alpha^{\ell} - \overline{\alpha}^{\ell}) = (\alpha^{\ell+1}-\overline{\alpha}^{\ell +1}) (\alpha^{k} - \overline{\alpha}^{k}).$$
It follows that
$$\alpha^{k} \overline{\alpha}^{\ell} (\overline{\alpha} - \alpha) = \alpha^{\ell} \overline{\alpha}^{k} (\overline{\alpha} - \alpha)$$
and therefore $(\alpha / \overline{\alpha})^{k - \ell} = 1$, but $\alpha / \overline{\alpha}$ is not a root of unity.
\end{proof}

\section{The three-dimensional case} \label{dimthree}

In this section we treat the method underlying Algorithm A. Below we prove the following result.

\begin{theorem} \label{constr}
Let $A = \{(x,y,z) \in \mathbb{Z}^3 : 0 \leq x < X, 0 \leq y < Y, 0 \leq z <Z\}$. Let $D$ be a finite set of nonproportional directions $(a_h,b_h,c_h)$ for $h=1,2,\dots, s$. Suppose the line sums of a function $f: A \to \mathbb{R}$ are given. Then a function $f^*: A \to \mathbb{R}$ can be constructed with the same line sums as $f$ has.
\end{theorem}

By the theory in \cite{ht07} this means that all such functions can be stated in parameterized form. In the next section we prove that $\mathcal{O}(sXYZ)$ elementary operations suffice to construct $f^*$.

Recall that $f$ is uniquely determined by its line sums in the directions of $D$ outside the union $T$ of its ghost domains.
Furthermore, if $(A,D)$ is valid, Equation \eqref{dimspace} provides the number of linearly independent functions $f^*: A \to \mathbb{R}$ with the same line sums as $f$ has and consequently such a number of function values of $f^*$ can be chosen. In the sequel the convex hull $C$ of $T$ plays an important role. Denote by $T^*$ the intersection of $T$ and the plane $z=Z-1$, and by $C^*$ the convex hull of $T^*$ in this plane.

\begin{prop} \label{hullz}
Suppose $(A,D)$ is valid. Then $T^*$ is a two-dimensional ghost domain obtained as a translation of the minimal ghost domain generated by the directions of the form $(a_h,b_h,0)$ in $D$. The translation is given by the vector $\sum_{c_h\neq0}\left(\abs{a_h}, \abs{b_h}, \abs{c_h}\right)$.
\end{prop}

\begin{proof}
By definition, $T^*$ consists of elements of $T$ with maximal $z$-value. In view of \eqref{elemt} and \eqref{ghostT} these elements satisfy $\varepsilon_h = 0$ if $c_h < 0$, and $\varepsilon_h = 1$ if $c_h > 0$. If $c_h = 0$ both $\varepsilon_h =0$ and $\varepsilon_h =1$ are possible, so that the boundary of $C^*$ contains line segments in the corresponding directions.
\end{proof}

\begin{cor} \label{numbers}
If $(A,D)$ is valid, then the minimal rectangular grid containing $C^*$ has size $X-\sum_{c_h \neq 0} ~|a_h|$ by $Y-\sum_{c_h \neq 0} ~|b_h|$. The number of points where the $f^*$-value can be freely chosen in the plane of $C^*$ equals $\left(X-\sum_h ~|a_h|\right) \left(Y-\sum_h ~|b_h|\right)$.
\end{cor}

\begin{proof}
For the first statement we recall that the directions $(a,b,c)$ with $c \not= 0$ do not contribute to $T^*$. The second statement follows from \eqref{numbers2}. It agrees with \eqref{dimspace} where the $z$-value has been fixed.
\end{proof}

\begin{figure}[htbp]
\centering
\subfigure[]
{\includegraphics[scale=.51]{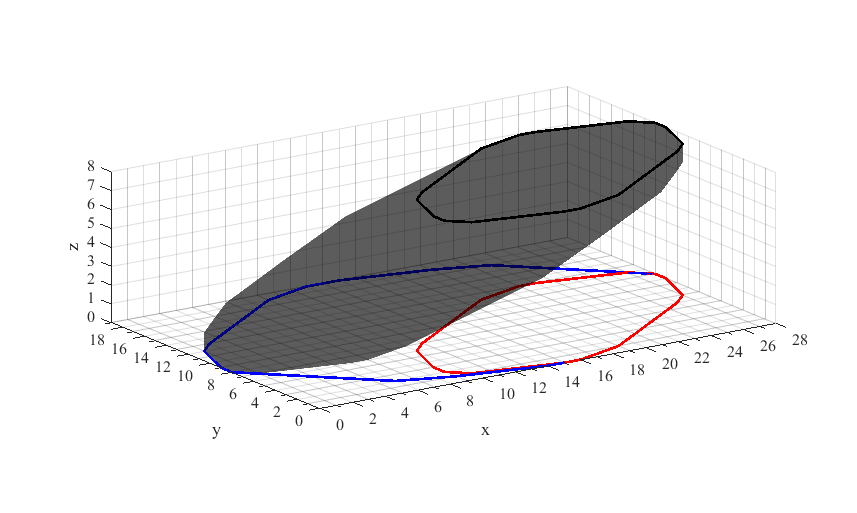}\label{fig_convhullxy_a}}
\hspace{.7cm}
\subfigure[]
{\includegraphics[scale=.6]{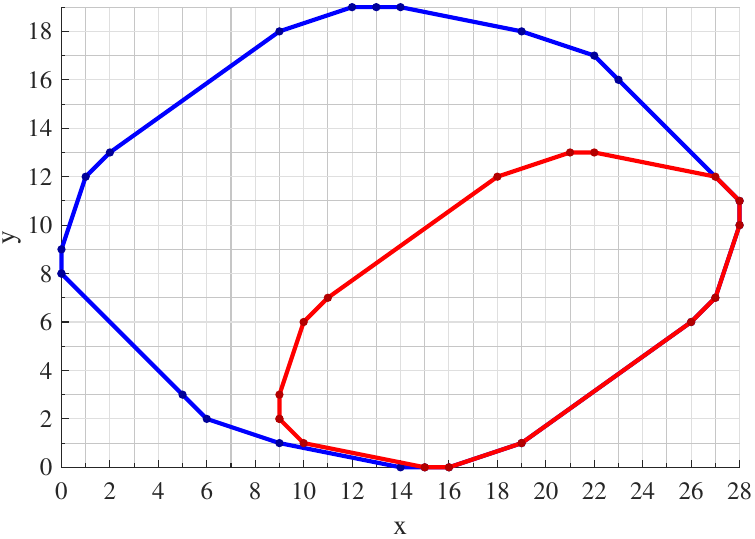}\label{fig_convhullxy_b}}
\caption{(Cf.~Example \ref{exg10}.) (a) The convex hull $C$ (shaded) of $\mathcal{G}_{10}$, the boundary of its projection on the $x$-$y$-plane along the $z$-axis (blue), the boundary of the convex hull $C^*$ (black) and its projection along the $z$-axis on the $x$-$y$-plane (red). (b) The boundary of the projection of $C$ along the $z$-axis (blue) and the boundary of the projection of $C^*$ along the $z$-axis (red).}
\end{figure}

\begin{ex} \label{exg10}
{\rm Consider the boundary ghost $\mathcal{G}_{10}$ generated by the directions $D = \{(0,1,0)$, $(1,0,0)$, $(1,1,0)$, $(1,0,1)$, $(-1,1,0)$, $(-3,-1,0)$, $(-3,1,-2)$, $(-1,-3,0)$, $(5,-1,0)$, $(5,-5,4)$, $(7,5,0)\}$. For $A$ we choose the grid $A := \left([0,29) \times [0,20) \times [0,9)\right) \cap \mathbb{Z}^3 $.
Since
\begin{displaymath}
\begin{array}{l}
\sum_h |a_h| = 28=X-1, \\
\sum_h |b_h| = 19=Y-1, \\
\sum_h |c_h| = 7=Z-2,
\end{array}
\end{displaymath}
we have a valid case. The number of free choices of $f^*$-values is 2 by \eqref{dimspace}. Suppose we know all line sums of a function $f:A\longrightarrow\mathbb{R}$ in the directions of $D$. Figure \ref{fig_convhullxy_a} shows $C$ (dark grey) and the boundary of $C^*$ (black). Figure \ref{fig_convhullxy_b} indicates their projections along the $z$-axis. The projection of $C$ is the blue curve and its interior, the projection of $C^*$ is the red curve and its interior.
In accordance with Corollary \ref{numbers} the size of the minimal grid containing $C^*$ is $X-\sum_{c_h \not= 0} ~|a_h|=20$ by $Y-\sum_{c_h \not= 0} ~|b_h|=14$ and at just $(29-28)(20-19)=1$ point of $A \cap \{z=Z-1\}$ the value of $f^*$ can be freely chosen.}
\end{ex}

\noindent {\it Proof of Theorem \ref{constr}.} The method uses Algorithm B several times and consists of four steps. If $(A,D)$ is nonvalid, then one can immediately proceed to Step 4. Steps 1-3 are used to reduce $A$ until the reduced $A$ is nonvalid for $D$. Of course, only $f^*$-values need to be computed which are not yet known. In the following steps, when projecting the convex hull $C$ of $T$ onto a coordinate plane, we omit writing the corresponding component in the projected direction.

\vspace{.3cm}
\noindent {\bf Step 1. Computation of $f$-values on lines with constant $x,z$.}
Suppose $(A,D)$ is valid. According to Corollary \ref{numbers}, $C^*$ is in a rectangular grid of size $X-\sum_{c_h \neq 0} ~|a_h|$ by $Y-\sum_{c_h \neq 0} ~|b_h|$. In Step 1 we compute the $f$-values on the $\sum_{c_h \neq 0} ~|a_h|$ lattice lines in $\{z=Z-1\}$ parallel to the $y$-axis which do not intersect $C^*$.

Consider the projection $C_y$ of $C$ along the $y$-axis onto the $x$-$z$-plane. We apply Algorithm B to the directions $D_y$ obtained from $D$ by omitting the second coordinates of the directions. The two-dimensional nonproportionality of $D_y$ follows from requirement 1) of Definition \ref{np}. By Theorem \ref{nonprop} we may apply Algorithm B with `primitive' replaced by `nonproportional'.

\begin{figure}
	\centering
		\includegraphics[scale=.7]{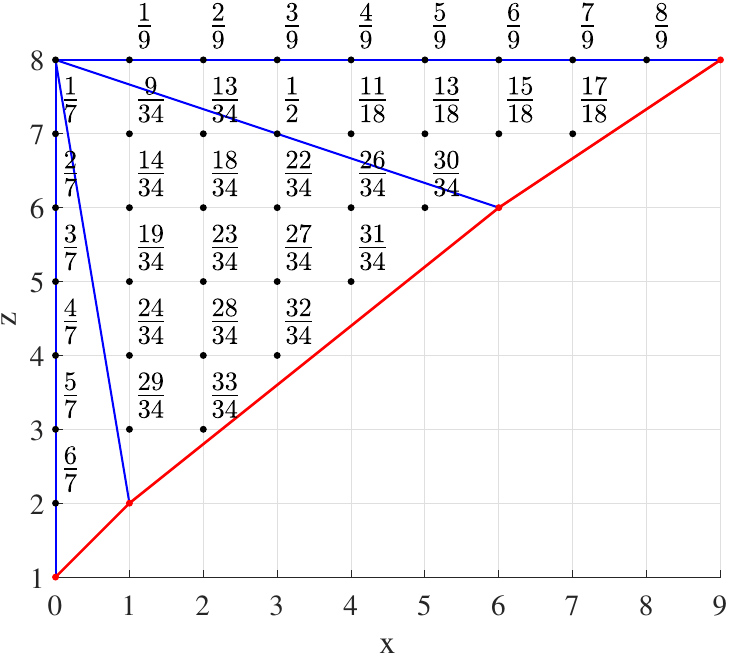}
	\caption{Projection of the convex hull $C$ onto the $x$-$z$-plane in the $y$-direction near the corner $(x,z) = (0,8)$. Weights are given for points outside the projection of $C$.}
	\label{fig_xzweights}
\end{figure}

As in Section \ref{sec:2dim} weights can be given to the points $\{(x,z) \in \mathbb{Z}^2 : 0 \leq x <X,  0 \leq z <Z\}$ outside $C_y$. These weights are such that, if in increasing ordering it is the turn of point $P_0=(x_0,z_0)$ located in the triangle along the side of the convex hull in the direction $(a_h,c_h)$ with $a_hc_h \not= 0$, then $P_0$ is the only point with unknown $f$-value on the lattice line through it in the direction $(a_h,c_h)$ (cf.~Figure \ref{fig_xzweights}). Therefore the lattice lines $(x_0,y,z_0) + t(a_h,b_h,c_h)~  (y \in \mathbb{Z}, t \in \mathbb{R})$ have no point in common with $C$. It follows that for every $y\in \mathbb{Z}$ the point $(x_0,y,z_0)$ is the only point on the lattice line $(x_0, y,z_0) + t(a_h,b_h,c_h),~t \in \mathbb{Z}$ with unknown $f$-value.
By the convexity of $C$ and Lemma \ref{lem:weight} the $f$-values of all the points $(x_0,y,z_0)$ with $0 \leq y < Y$ can be computed by subtracting from the line sum the $f$-values of all the other points contributing to the line sum. Thus the $f$-values of all the points $(x,y,z) \in A$ for which $(x,z)$ is outside $C_y$ can be computed. They cover the $\sum_{c_h \neq 0} ~|a_h|$ lattice lines in $\{z=Z-1\}$ parallel to the $y$-axis which do not intersect $C^*$.

\vskip.2cm
\noindent {\bf Example continued, Step 1.} As can be seen in Figure \ref{fig_convhullxy_b} the inner curve delimits $C^*$, which is contained in the rectangular grid $9 \leq x < 29, 0 \leq y < 14$. In Step 1 we compute the $f$-values of the points $(x,y,8)$ with $0 \leq x < 9, 0 \leq y < 20$. To do so, we consider the projection $C_y$ of $C$ in the $y$-direction. See Figure \ref{fig_xzweights}. The boundary of the projection of $C$ is indicated by the light coloured broken line from $(0,1)$ to $(9,8)$ consisting of segments in the direction (1,1), (5,4)  and (3,2). (These come from the directions $(1,0,1), (5,-5,4), (3,-1,2)$, i.e., the directions in $D$ with $ac \not= 0$.) The triangles and weights are as explained in Section \ref{sec:2dim}. The weights of $(x,z)$ are given to all points $(x,y,z)$ with $0 \leq y < 20$. Since the lines through these points in the directions $(1,0,1), (5,-5,4), (3,-1,2)$ do not intersect $C$, Algorithm B can be applied where, if it is the turn of point $(x_0,z_0)$, the $f$-values of all points $(x_0,y,z_0)$ with $0 \leq y < 20$ are computed. All these points are the only points on the corresponding lattice line with unknown $f$-value. Therefore this value can be found by subtraction. In this way we compute the $f$-values of the points $(x,y,8)$ with $0 \leq x < 9, 0 \leq y < 20$.

\begin{figure}
\centering
\includegraphics[scale=0.7]{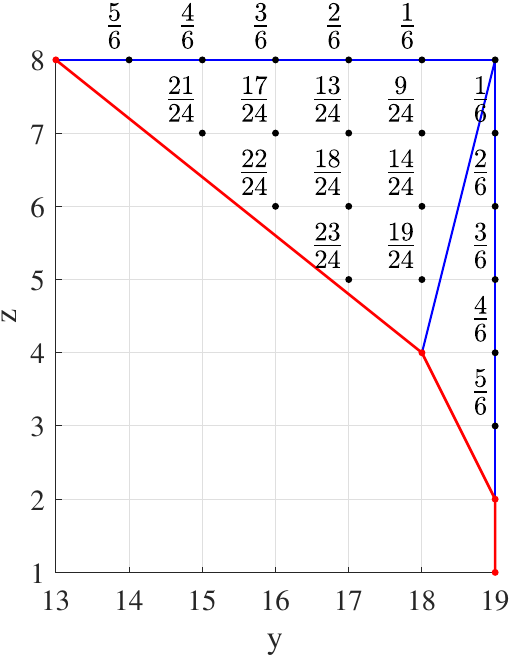}
\caption{Convex hull in the $y$-$z$-plane, with weights.}
\label{fig_yzweights}
\end{figure}

\vspace{.3cm}
\noindent {\bf Step 2. Computation of $f$-values on lines with constant $y,z$.}  We follow the same procedure as in Step 1 for the projection $C_x$ of $C$ in the $x$-direction. 
As a result, we retrieve all the $f$-values of the points $(x,y,z)\in A$ for which $(y,z)$ is outside $C_x$ (cf.~Figure \ref{fig_yzweights}). In particular we compute the $f$-values of all the points on the $\sum_{c_h \neq 0} ~|b_h|$ lattice lines in $\{z=Z-1\}$ parallel to the $x$-axis which do not intersect $C^*$.

\vskip.2cm
\noindent {\bf Example continued, Step 2}. In Step 2 we compute the $f$-values of the points $(x,y,8)$ with $9 \leq x < 29, 14 \leq y < 20$. To do so we consider the projection $C_x$ of $C$ in the $y$-direction. See Figure \ref{fig_yzweights}. The boundary of the projection of $C$ is indicated by the light coloured broken line from $(y,z) = (13,8)$ to $(19,2)$ consisting of segments in the directions $(5,-4)$ and $(1,-2)$. (These come from the directions $(5,-5,4), (3,-1,2)$, i.e., the directions in $D$ with $bc \not= 0$.) Again the extended Algorithm B can be applied where, if it is the turn of point $(y_0,z_0)$, the $f$-values of all points $(x,y_0,z_0)$ with $9 \leq x < 29$ are computed.

\begin{figure}
	\centering
		\includegraphics[scale=.7]{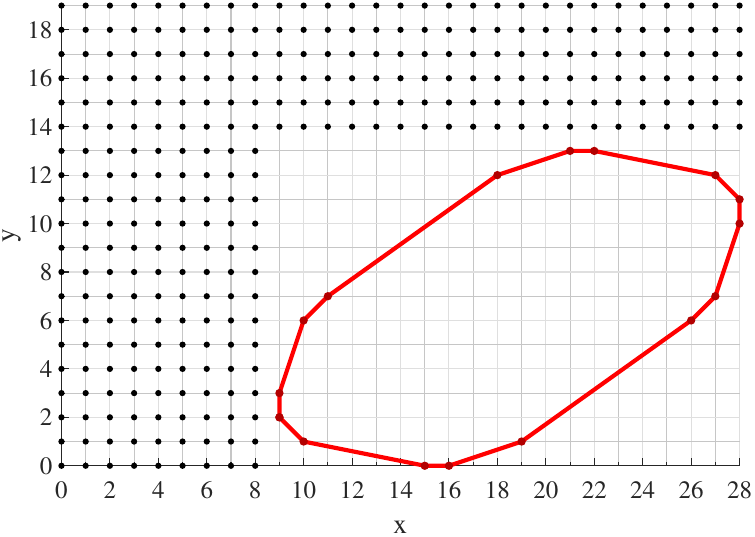}
	\caption{Intersection of the convex hull with the plane $z=8$. The isolated thick points have known $f$-values. In these points the $f^*$-value becomes the $f$-value.}
	\label{fig_reductionz5}
\end{figure}

\vspace{.3cm}
\noindent {\bf Step 3.} In this step we compute the $f^*$-values of the remaining points of the form $(x_0,y_0,Z-1)$. By Steps 1 and 2 for such a point there exist $x,y$ such that both $(x_0,y,Z-1)$ and $(x,y_0,Z-1)$ belong to $T$. Since the directions $(a_h,b_h,c_h) \in D$ with $c_h=0$ are nonproportional, we can apply the extended Algorithm B to compute the (remaining) $f^*$-values of all the points of the form $(x,y,Z-1) \in A$. In this process the number of $f^*$-values which can be freely chosen is given by Equation \eqref{numbers2} (cf.~Figure \ref{fig_reductionz5}). The choice of such values has to be performed first. Thereafter the remaining $f^*$-values can be computed. (Of course, every already found $f$-value is the $f^*$-value.)

We now remove the slice $A\cap\{z=Z-1\}$ from $A$ and consider the remaining part of $A$. The weights need not to be recalculated. If the new situation is valid, it follows from \eqref{elemt} and \eqref{ghostT} that the weights for $z=Z-2$ are as they were for $z=Z-1$. Again the number of $f^*$-values to be chosen is provided by Equation \eqref{numbers2}. After that, all remaining $f^*$-values of points with $z=Z-2$ can be calculated. This process is repeated for $z=Z-3$, $z=Z-4$, \ldots, $z=\sum_h |c_h|$. Equation \eqref{numbers2} gives each time the number of $f^*$-values which can be freely chosen. After completing this induction, $f^*$-values have been chosen in a number of points as Equation \eqref{dimspace} specifies. Different choices generate different solutions and the ultimately found solution depends only on these choices. What remains is a nonvalid rectangular grid $[0,X) \times [0,Y) \times [0, \sum_h |c_h|)$.

\vskip.2cm
\noindent {\bf Example continued, Step 3}.
For the third step, see Figure \ref{fig_reductionz5}. The isolated thick points are points with already known $f$-values. What is left is the rectangle $\{(x,y,8) : 9 \leq x <29, 0 \leq y < 14 \}$. We apply Algorithm B to the set of directions from $D$ with $z$-value 0, i.e., $(a,b)$ = $(0,1)$, $(1,0)$, $(1,1)$, $(-1,1)$ $(-3,-1)$, $(-1,-3)$, $(5,-1)$, $(7,5)$. Since $D$ is nonproportional, this set is two-dimensional nonproportional. The sum of the $|a|$ values is 19, the sum of the $|b|$ values is 13, so that, by \eqref{numbers2}, there is $(20-19)(14-13)=1$ choice for a $f^*$-value. We give the point $(9,2,8)$ an arbitrary $f^*$-value. Depending on this choice the $f^*$-values of the points $\{(x,y,8) : 9 \leq x <29, 0 \leq y < 14 \}$ can be computed using the extended Algorithm B.

Hereafter the same procedure is applied for the plane $z=7$. The changes are marginal: the scale for $z$ runs from 0 to 7 in place of from 0 to 8. This time the $f^*$-value of the point $(9,2,7)$ may be chosen in $\mathbb{R}$.

\begin{figure}
\centering
\includegraphics[scale=.8]{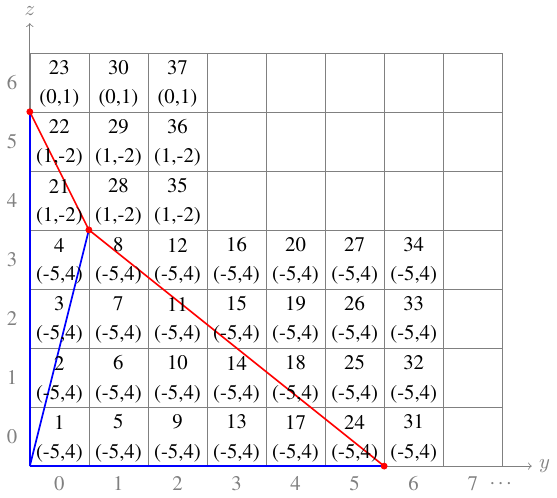}
\caption{Order (upper) and corresponding directions (lower) for the remaining $[0,29) \times [0,20) \times [0,7)$ block, projected along the $x$-axis. If it is the turn of a point $(y_0,z_0)$, then for all $x$ the points $(x,y_0,z_0)$ have a higher order number than all other points on the lattice line in the indicated direction. Thus the $f^*$-values in all points of $A'$ can be com\-put\-ed.}
\label{fig_yz_order}
\end{figure}

\vskip.2cm
\noindent {\bf Step 4.} We consider the remaining set $A'$ of points $(x,y,z)$ with $0 \leq x <X, 0 \leq y < Y, 0 \leq z < \sum_{h=1}^d |c_h|$. Then, $(A',D)$ is nonvalid. Since the directions $(a_h,b_h,c_h) \in D$ are nonproportional, the directions $(b_h,c_h)$ with $(a_h,b_h,c_h) \in D$ for some $a_h$ and $c_h \not= 0$ are two-dimensional nonproportional. We apply the extended Algorithm B to this set of directions. It provides an ordering of the points $(y,z) $ for which $(x,y,z) \in A'$ and for every such point a direction $(b_h,c_h)$ with $(a_h,b_h,c_h) \in D$ such that if it is the turn of $(y_0,z_0)$, then this point has higher weight than all other points on the lattice line $(y_0,z_0) + t(b_h,c_h)_{t \in \mathbb{Z}}$ (cf.~the proof of Theorem \ref{bROU} and Figure \ref{fig_yz_order}). Thus we can compute the $f^*$-values of all the points $(x,y_0,z_0)$ with $0 \leq x <X$ by considering the lattice line through $(x,y_0,z_0)$ in the direction $(a_h,b_h,c_h)$. After this step we know $f^*$.

\vskip.2cm
\noindent {\bf Example continued, Step 4}.
The fourth and last step is to deal with the $f^*$-values in the remaining block $[0,29) \times [0,20) \times [0,7)$. Since $\sum_n |c_n| =7$, we are now in the nonvalid case. We consider the projection of $C$ along the $x$-axis (see Figure \ref{fig_yz_order}). We apply Algorithm B to the directions $(0,1), (-5,4), (1,-2)$. The upper numbers in Figure \ref{fig_yz_order} show the ordering of the points $(y,z)$. After order number 20 there is periodicity modulo $7$. Following the points $(y,z)$ in increasing order, if it is the turn of $(y_0,z_0)$, the $f^*$-values can be computed for all points $(x,y_0,z_0)$ with $0 \leq x<29$, since $(x,y_0,z_0)$ is the only point without known $f^*$-value on the line in direction $(5,-5,4)$ if $0 \leq z_0 <4$, direction $(-3,1,-2)$ if $z_0=4$ or $5$ and direction $(1,0,1)$ if $z_0=6$. After this step we have found a function $f^*$ of which all the line sums in the directions of $D$ agree with those of $D$.

\section{Complexity}\label{sec:complexity}

Algorithm A deals with an $X\times Y\times Z$ grid and $s$ directions. As shown in \cite{cst}, Algorithm B involves $\mathcal{O}(sXY)$ elementary operations. Step 1 of the previous section requires $\mathcal{O}(sXZ)$ elementary operations to compute the weights and $\mathcal{O}(XYZ)$ operations to compute the $f^*$-values. Similarly, Step 2 requires $\mathcal{O}(sYZ + XYZ)$ elementary operations and the first part of Step 3 $\mathcal{O}(sXY + XYZ)$ elementary operations. Steps 1, 2 and 3 have to be repeated $Z - \sum_h |c_h|$ times. No new computation of the weights is required for Steps 1 and 2. No more than $XYZ$ $f^*$-values have to be calculated by a subtraction and an adjustment of the remaining line sum. Thus the reduction to a nonvalid situation takes in total $\mathcal{O}(sXYZ)$ elementary operations. Step 4 requires $\mathcal{O}(sYZ)$ elementary operations to compute the weights and $\mathcal{O}(XYZ)$ elementary operations to compute the $f^*$-values. All the order needed operations are of order at most $\mathcal{O}(sXYZ).$ We conclude that Algorithm A requires $\mathcal{O}(sXYZ)$ elementary operations.

\section{Conclusion}\label{sec:conclusion}

In this work we have presented a linear time approach for the reconstruction of a three-dimensional function $f: A \to \mathbb{R}$ from its line sums in a nonproportional set $D$ of directions. Here $\mathbb{R}$ may be replaced by any field or unique factorization domain. Such an approach was achieved through the extension to nonproportional directions of the two-dimensional algorithm presented in previous work and a multiple application of this extended algorithm. The definition of nonproportional is not symmetric in the coordinate directions. Of course, it is possible to interchange them to satisfy the definition.

It is desirable to have an algorithm like Algorithm A where the only restriction on the directions $(a,b,c)$ in $D$ is that they generate distinct one-dimensional vector spaces. One may define the three-dimensional weight of a point analogous to Lemma \ref{lem:weight}. In practice this is complicated, since it is awkward to characterize the convex hull $C$ in a useful way. Moreover, the line from the corner through the point need not intersect $C$. So new ideas are required.

\section*{Acknowledgments}
The first author's research has been supported from the Monash University Postgraduate Publications Award. The second author's research has been supported by D1 Research line of Università Cat\-to\-li\-ca del Sacro Cuore.

The authors thank the anonymous referee for his/her remarks, which have improved the paper considerably.

\bibliographystyle{plainurl}

\bibliography{cst_references_dim3_onlyneeded}

\end{document}